\nonstopmode
\def\cit#1#2{\ifx#1!\cite{#2}\else#2\fi}
\nonstopmode
\documentclass[]{amsart}
\usepackage{amssymb}
\usepackage{amscd}
\usepackage{graphicx}
\usepackage{psfrag}
\usepackage{verbatim}
\usepackage[all]{xy}
\newdir{ >}{{}*!/-10pt/@{>}}
\newdir^{ (}{{}*!/-5pt/@^{(}}
\newdir_{ (}{{}*!/-5pt/@_{(}}
\def\AMStag#1{\tag{#1}}
\def\undersetbrace#1\to#2{\underbrace{#2}_{#1}}
\def\oversetbrace#1\to#2{\overbrace{#2}^{#1}}
\def\AMSunderset#1\to#2{\underset{#1}{#2}}
\def\AMSoverset#1\to#2{\overset{#1}{#2}}
\def\therosteritem#1{{\rm (#1)}}
\def\East#1#2{\overset{#1}{\longrightarrow}}
\newenvironment{proclaim}[1]{\par\medskip\noindent{\bf #1.}\it}{\par\smallskip}
\newenvironment{demo}[1]{\par\smallskip{\bf #1.}}{\par\smallskip}
\def\nmb#1#2{\if0#1{#2}%
\else\if.#1{#2}%
\else{{#2}}%
\fi\fi}
\def\o{\circ}
\def\X{\mathfrak X}
\def\al{\alpha}

\def\ga{\gamma}

\def\ze{\zeta}

\def\vartheta{\theta}

\def\Ga{\Gamma}

\def\i{^{-1}}
\def\x{\times}

\def\Fl{\operatorname{Fl}}

\let\on=\operatorname
\def\Ad{\operatorname{Ad}}

\def\g{{\mathfrak g}}
\def\h{{\mathfrak h}}

\def\pr{{\operatorname{pr}}}

\def\p{\partial}
\def\today{\ifcase\month\or
 January\or February\or March\or April\or May\or June\or
 July\or August\or September\or October\or November\or December\fi
 \space\number\day, \number\year}
\def\AMSonly#1{}

\begin{document}
\title{Completing Lie algebra actions to Lie group actions
}
\author{Franz W.\  Kamber and
Peter W.\  Michor  }
\address{
Franz W. Kamber: 
Department of Mathematics, 
University of Illinois, 
1409 West Green Street, 
Urbana, IL 61801, USA}
\email{kamber@math.uiuc.edu}
\address{
P.\  W.\  Michor: Institut f\"ur Mathematik, Universit\"at Wien,
Strudlhofgasse 4, A-1090 Wien, Austria; {\it and:} 
Erwin Schr\"odinger Institut f\"ur Mathematische Physik,
Boltzmanngasse 9, A-1090 Wien, Austria}
\email{michor@esi.ac.at}
\thanks{
FWK and PWM were supported by `Fonds zur F\"orderung der wissenschaftlichen  
Forschung, Projekt P~14195~MAT'}

\begin{abstract} 
For a finite dimensional Lie algebra $\g$ of vector fields on a manifold $M$ we
show that $M$ can be completed to a $G$-space in a unversal way, 
which however is neither Hausdorff nor $T_1$ in general. Here $G$ is a 
connected Lie group with Lie-algebra $\g$. For a transitive $\g$-action the 
completion is of the form $G/H$ for a Lie subgroup $H$ which 
need not be closed. In general the completion can be constructed 
by completing each $\g$-orbit.
\end{abstract}
\def\LaTeXonly{}

\maketitle

\subsection*{\nmb0{1}. Introduction}
In \cit!{7}, Palais investigated when one could extend a local Lie group
action to a global one. 
He did this in the realm of non-Hausdorff
manifolds, since he showed, that completing a vector field $X$ on a Hausdorff
manifold $M$ may already lead to a non-Hausdorff manifold on which the
additive group $\mathbb R$ acts. We reproved this
result in \cit!{3}, being unaware of Palais' result. 
In \cit!{4} this result was extended to infinite dimensions and applied to
partial differential equations like Burgers' equation: Solutions of the PDE
were continued beyond the shocks and the universal
completion was identified. 

Here we give a detailed
description of the universal completion of a Hausdorff $\g$-manifold to a
$G$-manifold. For a homogeneous $\g$-manifold (where the finite dimensional
Lie algebra $\g$ acts infinitesimally transitive) we show that the
$G$-completion (for a Lie group $G$ with Lie algebra $\g$) is a homogeneous
space $G/H$ for a possibly non-closed Lie subgroup $H$ (theorem \nmb!{7}). 
In example \nmb!{8} we show that each such situation can indeed be
realized. For general $\g$-manifolds we show that one can complete each 
$\g$-orbit separately and replace the $\g$-orbits in $M$ by the resulting
$G$-orbits to obtain the universal completion ${_G}M$
(theorem \nmb!{9}). All $\g$-invariant structures on $M$ `extend' to $G$-invariant
structures on ${_G}M$.
The relation between our results and those of Palais are described in
\nmb!{10}. 

\subsection*{\nmb0{2}. $\g$-manifolds }
Let $\g$ be a Lie algebra. A $\g$-manifold is a (finite dimensional
Hausdorff) connected manifold $M$ together with a homomorphism of Lie
algebras $\ze=\ze^M:\g\to \X(M)$ into the Lie algebra of vector fields on
$M$. We may assume without loss that it is injective; if not replace $\g$
by $\g/\ker(\ze)$. We shall also say that $\g$ acts on $M$.

The image of $\ze$ spans an integrable distribution 
on $M$, which need not be of constant rank. So through each point of 
$M$ there is a unique maximal leaf of that distribution; we also call 
it the $\mathfrak g$-orbit through that point. It is an {\it initial 
submanifold} of $M$ in the sense that a mapping from a manifold into 
the orbit is smooth if and only if it is smooth into $M$, see 
\cit!{5},~2.14ff.

Let $\ell:G\x M\to M$ be a left action of a Lie group with Lie algebra
$\g$. Let $\ell_a:M\to M$ and $\ell^x:G\to
M$ be given by $\ell_a(x) = \ell^x(a) = \ell(a,x) = a.x$ for
$a\in G$ and $x\in M$.   
For $X\in \mathfrak g$ the 
{\it fundamental vector field} $\zeta_X=\ze_X^M\in \X(M)$ is given by 
$\ze_X(x) = -T_e(\ell^x).X = -T_{(e,x)}\ell.(X,0_x)=-\p_t|_0 \exp(tX).x$.
The minus sign is necessary so that $\ze:\g\to \X(M)$ becomes a Lie algebra
homomorphism.
For a right action the fundamental vector field mapping without minus 
would be a Lie algebra homomorphism. Since left actions are more
common, we stick to them.

\subsection*{\nmb0{3}. The graph of the pseudogroup }
Let $M$ be a $\g$-manifold, effective and connected, so that the 
action $\ze=\ze^M:\g\to \X(M)$ is injective.
Recall from \cit!{1},~2.3 that the pseudogroup $\Ga(\g)$ consists of 
all diffeomorphisms of the form
\begin{displaymath}
\Fl^{\ze_{X_n}}_{t_n}\o\dots
     \o\Fl^{\ze_{X_2}}_{t_2}\o\Fl^{\ze_{X_1}}_{t_1}|U
\end{displaymath}
where $X_i\in\g$, $t_i\in\mathbb R$, and $U\subset M$ are such that 
$\Fl^{\ze_{X_1}}_{t_1}$ is defined on $U$, $\Fl^{\ze_{X_2}}_{t_2}$ is 
defined on $\Fl^{\ze_{X_1}}_{t_1}(U)$, and so on.

Now we choose a connected Lie group $G$ with Lie algebra $\g$, 
and we consider the integrable distribution of constant rank 
$d =\dim (\g)$ on $G\x M$ which is given by 
\begin{equation}
\{(L_X(g),\ze^M_X(x)):(g,x)\in G\x M, X\in\g\}\subset TG\x TM,
\AMStag{\nmb.{3.1}}\end{equation}
where $L_X$ is the left invariant vector field on $G$ generated by 
$X\in\g$. This gives rise to the foliation $\mathcal F_\ze$ on $G\x M$, 
which we call the {\it graph foliation} of the $\g$-manifold $M$. 


Consider the following diagram, where $L(e,x)$ is the leaf through 
$(e,x)$ in $G\x M$, $\mathcal O_\g(x)$ is the $\g$-orbit through $x$ in 
$M$, and $W_x\subset G$ is the image of the leaf $L(e,x)$ in $G$. 
Note that $\pr_1 : L (e, x) \to W_x$ is a local diffeomorphism
for the smooth structure of $L (e, x).$ 
\begin{equation}
\xymatrix{
 & L(e,x)  \ar[dd] \ar[dr] \ar@{->>}[rr]^{\pr_2} &  & {\mathcal O}_\g (x) 
\ar@{^{(}->}[d] \\
 & & G\x M \ar[d]^{\pr_1} \ar[r]^{\pr_2} & M \\
[0,1] \ar@{-->}[uur]^{\tilde c} \ar[r]^{c}  &  W_x 
\ar@{^{(}->}[r]^{\text{open}} & G & }
\AMStag{\nmb.{3.2}}\end{equation}
Moreover we consider a piecewise smooth curve $c:[0,1]\to W_x$ with 
$c(0)=e$ and we assume that it is liftable to a smooth curve 
$\tilde c:[0,1]\to L(e,x)$ with $\tilde c(0)=(e,x)$. Its endpoint 
$\tilde c(1)\in L(e,x)$ does not depend on small (i.e.\ liftable to 
$L(e,x)$) homotopies of $c$ which respect the ends. This lifting 
depends smoothly on the choice of the initial point $x$ and gives 
rise to a local diffeomorphism 
$\ga_x(c): U\to \{e\}\x U\to \{c(1)\}\x U'\to U'$, 
a typical element of the pseudogroup $\Ga(\g)$ which is defined 
near $x$. 
See \cit!{1},~2.3 for more information and example \nmb!{4} below. 
Note, that the leaf $L(g,x)$ through $(g,x)$ is given by 
\begin{equation}
L(g,x)=\{(gh,y):(h,y)\in L(e,x)\}=(\mu_g\x\operatorname{Id})(L(e,x))
\AMStag{\nmb.{3.3}}\end{equation}
where $\mu:G\x G\to G$ is the multiplication and 
$\mu_g(h)=gh=\mu^h(g)$. 

\subsection*{\nmb0{4}. Examples }
It is helpful to keep the following examples in mind, which 
elaborate upon \cit!{1},~5.3. Let $G=\g=\mathbb R^2$, let $W$ be an 
annulus in $\mathbb R^2$ containing 0, and let $M_1$ be a simply 
connected piece of finite 
or infinite length of the universal cover of $W$. Then the Lie algebra 
$\g=\mathbb R^2$ acts on $M$ but not the group. Let $p:M_1\to W$ be the 
restriction of the covering map, a local diffeomorphism. 

Here $G \x_\g M_1 \cong G=\mathbb R^2$. Namely, the graph distribution is 
then also transversal to the fiber of $\pr_2:G\x M_1\to M_1$ (since the 
action is transitive and free on $M_1$), thus describes a principal 
$G$-connection on the bundle $\pr_2:G\x M_1\to M_1$. Each leaf is a covering 
of $M_1$ and hence diffeomorphic to $M_1$ since $M_1$ is simply connected. 
For $g\in \mathbb R^2$ consider 
$j_g:M_1 \East{{\operatorname{ins_g}}}{} \{g\}\x M_1 
     \subset G\x M_1 \East{{\pi}}{} G\x_\g M_1$
and two points $x\ne y\in M_1$.  We may choose a 
smooth curve $\ga$ in $M_1$ from $x$ to $y$, lift it into the leaf 
$L(g,x)$ and project it to a curve $c$ in $g+W$ from $g$ to 
$c(1)=g+p(y)-p(x)\in g+W$. Then $(g,x)$ and $(c(1),y)$ are on the same leaf.
So $j_g(x)=j_g(y)$ if and only if $p(x)=p(y)$. So we see that
$j_g(x)=g+p(x)$, and thus $G\x_\g M_1=\mathbb R^2$. This will also follow from
\nmb!{7}.

\centerline{
\psfrag{W}{$W$}
\psfrag{M1}{$M_1$}
\psfrag{M2}{$M_2$}
\includegraphics[width=7cm]{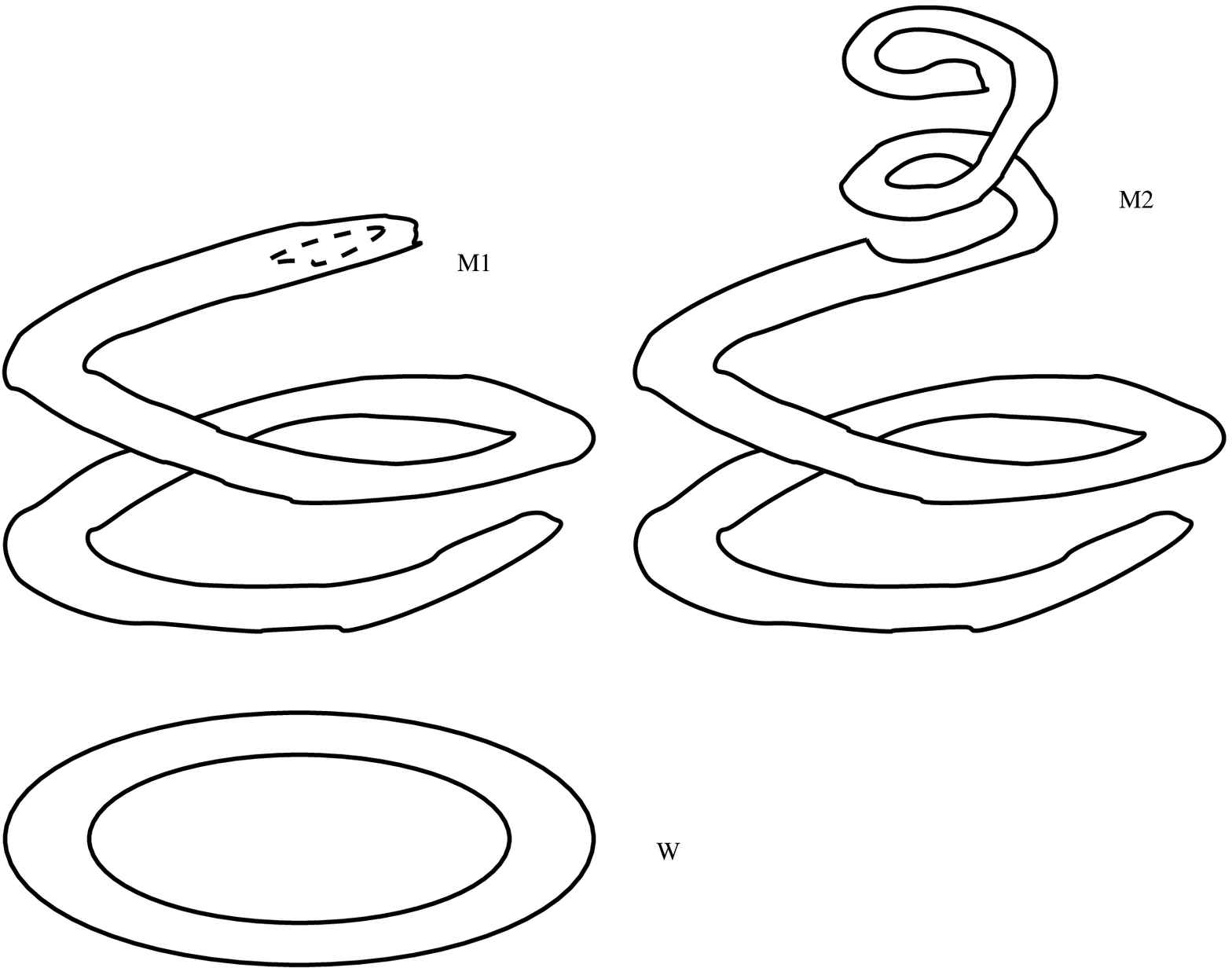}}
\medskip

Let us further complicate the situation by now omitting a small disk 
in $M_1$ so that it becomes non simply connected but still projects 
onto $W$, and let $M_2$ be a 
simply connected component of the universal cover of $M_1$ with the 
disk omitted. What happens now is that homotopic curves which act 
equally on $M_1$ act differently on $M_2$. 

It is easy to see with the methods described below that the completion
${_G}M_i=\mathbb R^2$ in both cases.  

\subsection*{\nmb0{5}. Enlarging to group actions }
In the situation of \nmb!{3} let us denote by 
$
{_G}M = G \x_\g M = G \times M / \mathcal F_{\zeta}
$ 
the space of leaves of the foliation $\mathcal F_{\zeta}$ on $G\x M$, 
with the quotient topology. For each $g\in G$ we consider the mapping 
\begin{equation}
j_g:M \East{{\operatorname{ins}_g}}{} \{g\}\x M 
     \subset G\x M \East{{\pi}}{} {_G}M  = G\x_\g M.
\AMStag{\nmb.{5.1}}\end{equation}
Note that the submanifolds $\{g\}\x M \subset G\x M$ are transversal 
to the graph foliation ${\mathcal F}_\ze.$ 
The leaf space ${_G}M$ of $G\x M$ of admits a unique smooth structure, 
possibly singular and non-Hausdorff, such that a mapping $f : {_G}M \to N$ 
into a smooth manifold $N$ is smooth if and only if the compositions 
$f \circ j_g : M \to N$ are smooth. For example we may use the 
structure of a {\it Fr\"olicher space} or {\it smooth space} induced by the
mappings $j_g$ in the
sense of \cit!{6},~section~23 on ${_G}M=G\x_\g M$.
The canonical open maps  $j_g : M \to {_G}M$ for $g \in G$ 
are called the charts of ${_G}M$~:
By construction, for each $x\in M$ and for $g'g\i$ near enough to $e$ in $G$ there
exists a curve $c:[0,1]\to W_x$ with $c(0)=e$ and $c(1)=g'g\i$ and
an open neighborhood $U$ of $x$ in $M$ such that for the smooth transformation
$\ga_x(c)$ in the pseudogroup $\Ga(\g)$ we have 
\begin{equation}
j_{g'}|U = j_{g}\o \ga_x(c).
\AMStag{\nmb.{5.2}}\end{equation}
Thus the mappings $j_g$ may serve as a replacement for charts in the
description of the smooth structure on ${_G}M$.
Note that the mappings $j_g$ are not injective in general. Even if $g=g'$ 
there might be liftable smooth loops $c$ in $W_x$ such that 
\thetag{\nmb|{5.2}} holds. Note also some similarity of the system of 
`charts' $j_g$ with the notion of an {\it orbifold} where one uses 
finite groups instead of pseudogroup transformations.

The leaf space ${_G}M = G \times_{\g} M$ is a smooth $G$--space where the
$G$-action is induced by 
$(g', x)\mapsto  (gg', x)$ in $G\x M$. 

\begin{proclaim}{Theorem}
The $G$--completion ${_G}M$ has the following {\it universal properties:}
\begin{enumerate}
\item[(\nmb.{5.3})]
  Given any Hausdorff $G$-manifold $N$ and $\g$-equivariant mapping 
  $f:M\to N$ there exists a unique $G$-equivariant continuous mapping 
  $\tilde f:{_G}M\to N$ with $\tilde f\o j_e = f$. Namely, the mapping 
  $\bar f:G\x M\to N$ given by $\bar f(g,x)=g.f(x)$ is smooth and 
  factors to $\tilde f:{}_GM\to N$.
\item[(\nmb.{5.4})]
  In the setting of \thetag{\nmb!{5.3}}, the universal property holds
  also for the $T_1$-quotient of ${_G}M,$ which is given as the quotient 
  $G\x M / \overline{\mathcal F}_\ze$ of $G\x M$ by the equivalence 
  relation generated by the closure of leaves. 
\item[(\nmb.{5.5})]
  If $M$ carries a symplectic or Poisson structure or a Riemannian metric 
  such that the
  $\g$--action preserves this structure or is even a Hamiltonian action
  then the structure `can be extended to ${}_GM$
  such that the enlarged $G$-action preserves these structures or is even
  Hamiltonian'.
\end{enumerate}
\end{proclaim}
 
\begin{demo}{Proof}
\thetag{\nmb!{5.3}} Consider the mapping 
$\bar f=\ell^N\o (\on{Id}_G\x f):G\x M\to N$ which is given by 
$\bar f(g,x)=g.f(x)$. Then by \thetag{\nmb!{3.1}} and \thetag{\nmb!{3.2}} we have for $X\in G$ 
\begin{align*}
T\bar f.(L_X(g),\ze^M_X(x)) 
&= 
T \ell.(L_X(g),T_x f.\ze^M_X(x)) 
\\&= 
T \ell.(R_{\Ad(g)X}(g),0_{f(x)}) + T\ell(0_g,\ze^N_X(f(x))) 
\\&= 
-\ze_{\Ad(g)X}(g.f(x)) + T\ell_g.\ze^N_X(f(x)) = 0.
\end{align*}
Thus $\bar f$ is constant on the leaves of the graph foliation on $G\x M$
and thus factors to $\tilde f:{_G}M\to N$. Since 
$\bar f(g.g_1,x)=g.g_1.f(x)=g.\bar f(g,x)$, the mapping $\tilde f$ is
$G$-equivariant. Since $N$ is Hausdorff, $\tilde f$ is even constant on the
closure of each leaf, thus 
\thetag{\nmb!{5.4}} holds also.

\thetag{\nmb!{5.5}} Let us treat Poisson structure $P$ on $M$. For
symplectic structures or Riemannian metrics the argument is similar and simpler.
Since the Lie derivative along fundamental vector
fields of $P$ vanishes, the pseudogroup transformation $\ga_x(c)$ in
\thetag{\nmb!{5.2}} preserves $P$. Since ${}_GM$ is the quotient of the
disjoint union of all spaces $\{g\}\x M$ for $g\in G$ under the equivalence
relation described by \thetag{\nmb!{5.2}}, $P$ `passes down to this quotient'.
Note that we refrain from putting too much meaning on this statement.
\qed\end{demo}

The universal property \thetag{\nmb!{5.3}} holds also for
smooth $G$-spaces $N$ which need not be Hausdorff, nor $T_1$, but should
have tangent spaces and foliations so that it is meaningful to talk about
$\g$-equivariant mappings. We will not go into this, but see \cit!{6},
section 23 for some concepts which point in this direction.

As an application of the universal property of the $G$--completion 
${_G}M,$ we see that ${_G}M$ depends on the choice of $G$ in the 
following way. We write $G = \Ga \backslash \widetilde{G},$ where 
$\widetilde{G}$ is the simply connected Lie group with Lie algebra $\g$ 
and $\Ga \subset \widetilde{G}$ is the discrete central subgroup such that 
$\Ga \cong \pi_1 (G).$ Then we have ${_G}M \cong \Ga \backslash 
{_{\widetilde G}}M$ as $G$--spaces, so that ${_{\widetilde G}}M$ 
is potentially less singular than ${_G}M.$ 

\subsection*{\nmb0{6}. Example}
Let $\g=\mathbb R^2$ with basis $X,Y$, let 
$M=\mathbb R^3 \setminus \{(0,0,z):z \in \mathbb R\}$, and let 
$\zeta^\al : \g \to \X(M)$ be given by 
\begin{equation}
\zeta^\al_X=\partial_{x}+\al ~\frac{y z}{x^2+y^2} ~\partial_{z}
\qquad, \qquad 
\zeta^\al_Y=\partial_{y}-\al ~\frac{x z}{x^2+y^2} ~\partial_{z}~, 
~\al > 0
\AMStag{\nmb.{6.1}}\end{equation}
which satisfy $[\zeta^\al_X,\zeta^\al_Y] = 0$. 
By construction of the graph foliation ${\mathcal F}_{\ze^\al}$ in 
\thetag{\nmb!{3.1}} and the procedure summarized in diagram \thetag{\nmb!{3.2}}, 
the leaves of ${\mathcal F}_{\ze^\al}$ are determined explicitly as follows.  
For any smooth curve $c (t) = (\xi (t), \eta (t)) \in G$ starting at 
$(\xi_0, \eta_0)$ we have 
$\dot{c} (t) = \dot{\xi} (t)~X + \dot{\eta} (t)~Y \in \g$ and the lifted 
curve $(c (t), \mathbf{y} (t))$ is in the leaf $L ((\xi_0,\eta_0), \mathbf{y}_0)$ 
if and only if it satisfies the first order ODE 
\begin{equation}
(\mathbf{y} (t), \dot{\mathbf{y}} (t)) = 
\dot{\xi} (t)~\zeta^\al_X (\mathbf{y} (t)) + 
\dot{\eta} (t)~\zeta^\al_Y (\mathbf{y} (t))
\AMStag{\nmb.{6.2}}\end{equation}
with initial value $\mathbf{y} (0) = \mathbf{y}_0 = (x_0, y_0, u=z_0) \in M$. 
Substituting \thetag{\nmb!{6.1}} into \thetag{\nmb!{6.2}}, we see that this 
ODE is linear, that is $\dot{x} = \dot{\xi}~, ~\dot{y} = \dot{\eta}$ 
and $\dot{z} = - \al ~z ~\frac{x \dot{\eta} - y \dot{\xi}}{r^2} = 
- \al ~z ~\frac{x \dot{y} - y \dot{x}}{r^2}$, where $r^2 = x^2+y^2.$ 
Thus the projection $\mathbf{x} (t)$ of $\mathbf{y} (t)$ to the 
$(x,y)$--plane is given by $\mathbf{x} (t) = c (t) - ((\xi_0,\eta_0) - \mathbf{x}_0) 
= c (t) - (\xi_0 - x_0 , \eta_0 - y_0),$ whereas the third equation leads to 
\begin{equation}
z (t) = u ~e^{-\al \int_0^t d\theta} = u ~e^{-\al (\theta (t) - \theta_0)} 
= u e^{\al \theta_0} ~e^{-\al \theta (t)}~, 
\AMStag{\nmb.{6.3}}\end{equation}
where $\theta$ is the angle function in the $(x,y)$--plane. 
This depends only on the endpoints $\mathbf{x}_0~, ~\mathbf{x} (t)$ and 
the winding number of the curve $\mathbf{x}$ and is otherwise 
independent of $\mathbf{x}.$ 
Incompleteness occurs whenever the curve $\mathbf{x}$ goes to 
$(0,0) \in {\mathbb R}^2$ in finite time $\bar{t}<\infty,$ that is  
$\mathbf{x} (t) \to (0,0)~, ~t \uparrow \bar{t}$ or equivalently 
$c (t) \to (\xi_0,\eta_0) - \mathbf{x}_0~, ~t \uparrow \bar{t}.$ 
It follows that the leaf $L ((\xi_0,\eta_0), \mathbf{y}_0)$
is parametrized by $(r, \theta)\in {\mathbb R}_+\x {\mathbb R}$ 
with $z=z(\theta)$ being independent of $r>0$ and that 
\begin{equation}
\pr_1:L ((\xi_0,\eta_0), \mathbf{y}_0) \to W_{(\xi_0,\eta_0),\mathbf{y}_0} = 
{\mathbb R}^2 \setminus \{(\xi_0,\eta_0) - \mathbf{x}_0 \}
\AMStag{\nmb.{6.4}}\end{equation}
in \thetag{\nmb!{3.2}} is a universal covering. 
This is visibly consistent with \thetag{\nmb!{3.3}}. 
In order to pa\-ra\-met\-rize the space of leaves ${_G}M$, 
we observe that the parameter $\mathbf{x}_0$ can be eliminated. 
In fact, from the previous formulas we see that  
\begin{equation}
L ((\xi_0',\eta_0'), (\mathbf{x}_0', u')) = 
L ((\xi_0,\eta_0), (\mathbf{x}_0, u)), 
\AMStag{\nmb.{6.5}}\end{equation}
if and only if $(\xi_0',\eta_0') - \mathbf{x}_0' = 
(\xi_0,\eta_0) - \mathbf{x}_0$ and $u' = u e^{\al (\theta_0 -\theta_0')},$ 
so that we have 
$z' (\theta) = u' e^{\al \theta_0'} ~e^{-\al \theta (t)} = 
u e^{\al \theta_0} ~e^{-\al \theta (t)} = z (\theta).$
In particular, it follows that 
\begin{equation}
L ((\xi_0,\eta_0), \mathbf{y}_0) = L ((\xi_0'+1,\eta_0'), (1, 0, u')), 
\AMStag{\nmb.{6.6}}\end{equation}
where $(\xi_0',\eta_0') = (\xi_0, \eta_0) - \mathbf{x}_0~, ~
u' = ue^{\al \theta_0}~, ~\theta_0'=0,$ projecting to 
${\mathbb R}^2 \setminus \{(\xi_0',\eta_0') \}.$  
Therefore the leaves of the form 
$L ((\xi_0+1,\eta_0), (1, 0, u))$ 
are distinct for different values of $(\xi_0,\eta_0)$ and fixed value of 
$u$ and from the relation \thetag{\nmb!{3.3}} we conclude that 
\begin{equation}
L ((\xi_0+1,\eta_0), (1, 0, u))= (\xi_0,\eta_0) + L ((1, 0), (1, 0, u)),
\AMStag{\nmb.{6.7}}\end{equation}
that is $G = {\mathbb R}^2$ acts without isotropy on ${_G}M$. 
We also need to determine the range for the parameter $u$. 
Obviously, we have $L ((1, 0), (1, 0, u')) = 
L ((1, 0), (1, 0, u))$ if and only if 
$u' = e^{2 \pi \al n} u$ for $n \in \mathbb Z$. 
Thus these leaves are parametrized by 
$[u],$ taking values in the quotient of the 
additive group $\mathbb R$ under the multiplicative group 
$\{ e^{2 \pi \al n }: n \in \mathbb Z \},$ that is 
\begin{equation}
\{ 0 \}\cup {\mathbb S}^1_{+} \cup {\mathbb S}^1_{-} \cong 
\{ 0 \} \cup {\mathbb R}^\x_+ / \{ e^{2 \pi \al n}:n \in \mathbb Z\} 
\cup {\mathbb R}^\x_- / \{ e^{2 \pi \al n}:n \in \mathbb Z\}.
\AMStag{\nmb.{6.8}}\end{equation}
The topology on the above space is determined by the leaf closures, 
respectively the orbit closures. First we have 
$
\overline{L ((\xi_0+1,\eta_0), (1, 0, u))} = 
(\xi_0,\eta_0) + \overline{L ((1, 0), (1, 0, u))}
$
in $G \x M$ and it is sufficient to determine the 
closures of $L ((1, 0), (1, 0, u)).$ 
For $(1,0,u)\in M$ with $u \ne 0$ we consider the curve 
$c(\theta)= e^{i \theta}\in G=\mathbb R^2.$ 
It is liftable to $G \x M$ and determines on $M$ the curve 
$\mathbf{y}(t)= (\cos \theta, \sin \theta, u e^{-\al \theta}).$ 
Thus the curve $(c(\theta), \mathbf{y}(\theta))$ in the leaf through 
$(1,0;1,0,u) \in G\x M\subset \mathbb R^5$ has a limit cycle for 
$\theta \to \infty$ which lies in the different leaf through 
$(1,0;1,0,0)$ which is closed, given by the $(x,y)$--plane 
$(\mathbb R^2\x 0) \setminus 0$ at level $(1,0) \in G.$
Thus we have 
\begin{equation}
\overline{L ((1, 0), (1, 0, u))} = 
L ((1, 0), (1, 0, u)) \cup L ((1, 0), (1, 0, 0)).
\AMStag{\nmb.{6.9}}\end{equation}
Hence the leaf $L ((1, 0), (1, 0, u))$ is not closed and the 
topological space ${_G}M$ is not $T_1$ and not a manifold. 
The orbits of the $\g$-action are determined by the leaf structure 
via $\pr_2$ in diagram \thetag{\nmb!{3.2}} and they look here as follows: 
The $(x,y)$--plane $(\mathbb R^2\x 0)\setminus 0$ is a closed orbit. 
Orbits above this plane are helicoidal staircases leading down 
and accumulating exponentially at the $(x,y)$--plane. 
Orbits below this plane are helicoidal staircases leading up and 
again accumulating exponentially.
Thus the orbit space $M/\g$ of the $\g$-action is given by
\thetag{\nmb!{6.8}}, with the point $0$ being closed. 
By \thetag{\nmb!{6.9}}, the closure of any orbit represented by a 
point $[u]$ on one of the circles is given by $\{ [u], 0 \}$. 
>From \thetag{\nmb!{6.6}} and \thetag{\nmb!{6.7}}, 
we see that the $G$-completion ${_G}M$ has a section over 
the orbit space ${_G}M / G \cong M / \g$ given by 
$[u] \mapsto L ((1, 0), (1, 0, u)).$ 
Therefore ${_G}M \cong G \x M / \g = {\mathbb R}^2 \x 
\{ ~\{ 0 \}\cup {\mathbb S}^1_{+} \cup {\mathbb S}^1_{-}~\}.$  

The structure of the completion and the orbit spaces are independent of 
the deformation parameter $\al>0$ in \thetag{\nmb!{6.1}}. 
However for $\al \downarrow 0,$ the completion just means adding 
in the $z$--axis, that is we get ${_G}M \cong {\mathbb R}^3$ with 
$G={\mathbb R}^2$ acting by parallel translation on the affine planes 
$z=c,$ and $M / \g \cong {_G}M / G \cong {\mathbb R}$ as it should be. 

It was pointed out to us \cit!{2} that one can make this example still more
pathological:  Consider the above example only in a cylinder over the
anulus $0<x^2+y^2<1$. 
Add an open handle to the disk and continue the
$\mathbb R^2$-action on the cylinder over the disk with an open handle added 
in such a way that there is a shift in the
$z$--direction when one traverses the handle. Then one of the
helicoidal staircases is connected to the the disk itself, so it
accumulates onto itself. This is called a `resilient leaf' in foliation
theory. 

\begin{proclaim}{\nmb0{7}. Theorem}
Let $M$ be a connected transitive effective $\g$-manifold. Let $G$ be a
connected Lie group with Lie algebra $\g$. 
Then we have:
\begin{enumerate}
\item[(\nmb.{7.1})] 
  Then there exists a subgroup
  $H\subset G$ such that the $G$-completion ${_G}M$ is diffeomorphic to
  $G/H$.
\item[(\nmb.{7.2})] 
  The Hausdorff quotient of ${_G}M$ is the homogeneous manifold 
  $G/\overline H$. It has the following universal property: 
  For each smooth $\g$-equivariant mapping $f:M\to N$ into a Hausdorff
  $G$-manifold $N$ there exists a unique smooth $G$-equivariant mapping
  $\tilde f: G/\overline H\to N$ with 
  $f=\tilde f\o \pi\o j_e:M\to G/H \East{{\pi}}{} G/\overline H\to N$. 
\item[(\nmb.{7.3})] 
  For each leaf $L(g,x_0)\subset G\x M$ the projection 
  $\pr_2:L(g,x_0)\to M$ is a smooth fiber bundle with typical fiber $H$.
\end{enumerate}
\end{proclaim}

\begin{demo}{Proof}
Since the action is transitive we have the exact sequence of vector bundles 
over $M$ 
\begin{displaymath}
0\to \on{iso} \to M\x \g \East{{\hat\ze}}{} TM \to 0.
\end{displaymath}

\thetag{\nmb!{7.1}}
We choose a base point $x_0\in M$. 
The $G$-completion is given by ${_G}M = G\x_\g M$, the orbit space of the
$\g$-action on $G\x M$ which is given by $\g\ni X\mapsto L_X\x \ze^M_X$,
and the $G$-action on the completion is given by multiplication from the
left. 
The submanifold $G\x\{x_0\}$ meets each $\g$-orbit in $G\x M$ transversely,
since  
\begin{align*}
T_{(g,x_0)}(G\x \{x_0\}) + T_{(g,x_0)}L(g,x_0) 
&= \{L_X(g)\x 0_{x_0} + L_Y(g)\x \ze_Y(x_0): X,Y\in \g\}
\\&=
T_{(g,x_0)}(G\x M).
\end{align*}
By \thetag{\nmb!{3.3}} we have $L(g,x)=g.L(e,x)$ so that we can define the linear
subspace $\g_{x_0}=\h\subset \g$ by 
\begin{align*}
X\in \h\quad \Longleftrightarrow \quad &
X \x 0_{x_0}\in T_{(e,x_0)}(G\x \{x_0\}) \cap T_{(e,x_0)}L(e,x_0) 
\\\Longleftrightarrow \quad &
L_X(g)\x 0_{x_0}\in T_{(g,x_0)}(G\x \{x_0\}) \cap T_{(g,x_0)}L(g,x_0) 
\end{align*}
Since $G\x\{x_0\}$ is a leaf of a foliation and the $L(e,x)$ also form a
foliation, $\h$ is a Lie subalgebra of $\g$. Let $H_0$ be the connected Lie
subgroup of $G$ which corresponds to $\h$. Then clearly 
$H_0\x \{x_0\}\subset G\x \{x_0\}\cap L(e,x_0)$. 
Let the subgroup $H\subset G$ be given by 
\begin{displaymath}
H=\{g\in G: (g,x_0)\in L(e,x_0)\}
=\{g\in G: L(g,x_0)= L(e,x_0)\},
\end{displaymath}
then the $C^\infty$-curve component of $H$ containing $e$ is just $H_0$.
So $H$ consists of at most countably many $H_0$-cosets. Thus $H$ is a Lie
subgroup of $G$ (with a finer topology, perhaps).
By construction the orbit space $G\x_{\g}M$ equals the quotient of the
transversal $G\x \{x_0\}$ by the relation induced by intersecting with leaves
$L(g,x_0)$, i.e., $G\x_\g M= G/H$.

\thetag{\nmb!{7.2}}
Obviously the $T_1$-quotient of $G/H$ equals the Hausdorff quotient
$G/\overline H$ which is a smooth manifold. The universal property is easily
seen. 

\thetag{\nmb!{7.3}}
Let $x\in M$ and $(g,x)\in L(e,x_0)=L(g,x)=g.L(e,x)$. So it suffices to
treat the leaf $L(e,x)$. 
We choose $X_1,\dots X_n\in \g$ such that
$\ze_{X_1}(x),\dots,\ze_{X_n}(x)$ form a basis of the tangent space $T_xM$. 
Let $u:U\to \mathbb R^n$ be a chart on $M$ centered at $x$ such that $u(U)$ is
an open ball in $\mathbb R^n$ and such that $\ze_{X_1}(y),\dots,\ze_{X_n}(y)$
are still linearly independent for all $y\in U$. For $y\in U$ consider the
smooth curve $c_y:[0,1]\to U$ given by $c_y(t)=u\i(t.u(y))$. We
consider 
\begin{align*}
\p_t c_y(t)&=c_y'(t)={\textstyle\sum}_{i=1}^n f_y^i(t)\,\ze_{X_i}(c_y(t)),\quad
  f_y^i\in C^\infty([0,1],\mathbb R)\\
X_y(t) &= {\textstyle \sum}_{i=1}^n f_y^i(t)\,X_i\in \g,\qquad
  X\in C^\infty([0,1],\g) \\
g_y &\in C^\infty([0,1],G),\quad T(\mu_{g_y(t)})\p_t g_y(t) = X_y(t),\quad
  g_y(0)=e,
\end{align*}
and everything is also smooth in $y\in U$.
Then for $h\in H$ we have $(h.g_y(t),c_y(t))\in L(e,x)$ since 
\begin{displaymath}
\p_t (h.g_y(t),c_y(t)) = (L_{X_y(t)}(h.g_y(t)), \ze_{X_y(t)}(c_y(t))).
\end{displaymath}
Thus $U\x H\ni (y,h)\mapsto \pr_2\i(U)\cap L(e,x)$ is the required fiber
bundle parameterization. 
\qed\end{demo}

\subsection*{\nmb0{8}. Example}
Let $G$ be simply connected Lie group and let $H$ be a connected Lie group
of $G$ which is not closed. For example, let $G=Spin(5)$ which is compact of 
rank 2 and let $H$ be a dense 1-parameter subgroup in its 2-dimensional
maximal torus. Let $\on{Lie}(G)=\g$ and $\on{Lie}(H)=\h$. 
We consider the foliation of $G$ into right $H$-cosets $gH$ which is
generated by $\{L_X: X\in \h\}$ and is left invariant under $G$.
Let $U$ be a chart centered at $e$ on $G$ which is adapted to this
foliation, i.e. $u:U\to u(U)= V_1\x V_2\subset \mathbb R^k\x\mathbb R^{n-k}$ such
that the sets $u\i(V_1\x\{x\})$ are the leaves intersected with $U$. 
We assume that $V_1$ and $V_2$ are open balls, and that $U$ is so small
that $\exp: W\to U$ is a diffeomorphism for a suitable convex open set
$W\subset \g$. Of
course $\g$ acts on $U$ and respects the foliation, so this $\g$-action
descends to the leave space $M$ of the foliation on $U$ which is
diffeomorphic to $V_2$. 

\begin{proclaim}{Lemma} In this situation, for the $G$-completion we
have $G\x_\g M=G/H$  
\end{proclaim}

\begin{demo}{Proof}
We use the method described in the end of the proof of theorem \nmb!{7}:
${_G}M=G\x_\g M$ is the quotient of the transversal $G\x\{x_0\}$ by the
relation induced by intersecting with leaves $L(g,x_0)$. Thus we have
to determine the subgroup $H_1=\{g\in G:(g,x_0)\in L(e,g)\}$. 

Obviously any smooth curve $c_1:[0,1]\to H$ starting at $e$ is liftable to
$L(e,x_0)$ since it does not move $x_0\in M$. So $H\subseteq H_1$, and
moreover $H$ is the $C^\infty$-path component of the identity in $H_1$. 

Conversely, if $c=(c_1,c_2):[0,1]\to L(e,x_0)\subset G\x M$ is a smooth curve 
from $(e,x_0)$ to $(g,x_0)$ then $c_2$ is a smooth loop through $x_0$ in
$M$ and there exists a smooth homotopy $h$ in $M$ which contracts $c_2$ to
$x_0$, fixing the ends. Since $\pr_2:L(e,x_0)\to M$ is a fiber bundle by
\thetag{\nmb!{7.3}} we
can lift the homotopy $h$ from $M$ to $L(e,x_0)$ with starting curve $c$,
fixing the ends, and deforming $c$ to a curve $c'$ in $L(e,x_0)\cap
\pr_2\i(x_0)$. Then $\pr_1\o c'$ is a smooth curve in $H_1$ connecting $e$
and $g$. 

Thus $H_1=H$, and consequently ${_G}M=G/H$.
\qed\end{demo}

\begin{proclaim}{\nmb0{9}.  Theorem}
Let $M$ be a connected $\g$-manifold. Let $G$ be a connected Lie group with Lie
algebra $\g$. Then the $G$-completion ${_G}M$ can be described in the
following way: 
\begin{enumerate}
\item[(\nmb.{9.1})] 
  Form the leaf space $M/\g$, a quotient of $M$ which may be
  non-Hausdorff and not $T_1$ etc.
\item[(\nmb.{9.2})] 
  For each point $z\in M/\g$, replace the orbit $\pi\i(z)\subset M$ by
  the homogeneous space $G/H_x$ described in theorem \nmb!{7}, where $x$
  is some point in the orbit $\pi\i(z)\subset M$. One can use transversals
  to the $\g$-orbits in $M$ to describe this in more detail. 
\item[(\nmb.{9.3})] 
  For each point $z\in M/\g$, one can also replace the orbit 
  $\pi\i(z)\subset M$ by the homogeneous space $G/\overline{H_x}$ described 
  in theorem \nmb!{7}, where $x$ is some point in the orbit $\pi\i(z)\subset M$. 
  The resulting $G$-space has then Hausdorff orbits which are smooth
  manifolds, but the same orbit space as $M/\g$. 
\end{enumerate}
\end{proclaim}

See example \nmb!{6} above.

\begin{demo}{Proof}
Let $\mathcal O(x)\subset M$ be the $\g$-orbit through $x$, i.e., the leaf
through $x$ of the singular foliation (with non-constant leaf dimension) on $M$ 
which is induced by the $\g$-action. 
Then the $G$-completion of the orbit $\mathcal O (x)$ is ${_G}\mathcal O(x)=G/H_x$ 
for the Lie subgroup $H_x\subset G$ described in theorem
\thetag{\nmb!{7.1}}. By the universal property of the $G$-completion we get a
$G$-equivariant mapping ${_G}\mathcal O(x)\to {_G}M$ which is injective and a
homeomorphism onto its image, since we can repeat the construction of
theorem \thetag{\nmb!{7.1}} on $M$. 
Clearly the mapping $j_e:M\to {_G}M$ induces a homeomorphism between the orbit
spaces $M/\g\to {_G}M/G$.

Now let $s: V\to M$ be an embedding of a submanifold which is a transversal
to the $\g$-foliation at $s(v_0)$: We have $Ts\cdot T_{v_0}V\oplus
\ze_{s(v_0)}(\g)=T_{s(v_0)}M$. Then $s$ induces a mapping $V\to G\x M$ and
$V\to {_G}M$ and we may use the point $s(v)$ in replacing $\mathcal O(s(v))$ by
$G/H_{s(v)}$ for $v$ near $v_0$.
\qed\end{demo}

\medbreak
The following diagram summarizes the relation between the preceding 
constructions. 
\begin{equation}
\xymatrix{
M \ar[rr] \ar@{<->}[d]^{=} & & \bigcup_{[x] \in M/\g} ~G / H_x 
\ar[d]^{\cong} \ar@{->>}[r] &
\bigcup_{[x] \in M/\g} ~G/\overline{H_x} \ar@{-->>}[d] \ar@{-->>}[ddl] \\
M \ar@{->>}[d]^{\pi} \ar[rr]^{j_e} & & {_G}M = G \x_\g M 
\ar@{->>}[d]^{\pi_G}  \ar@{->>}[r]  & G\x M / \overline{\mathcal F}_\ze 
\ar@{->>}[d]^{\bar{\pi}_G} \\
M / \g \ar[rr]^{\cong} & &  {_G}M / G \ar@{->>}[r] & 
(G\x M / \overline{\mathcal F}_\ze)/G }
\AMStag{\nmb.{9.4}}\end{equation}
Note that taking the $T_1$--quotient $G\x M / \overline{\mathcal F}_\ze$ 
of the leaf space ${_G}M$ may be a very severe reduction. 
In example \nmb!{6} the isotropy groups $H_x$ 
are trivial and we have $G\x M / \overline{\mathcal F}_\ze = {\mathbb R}^2 \x \{ 0\}$ 
and $(G\x M / \overline{\mathcal F}_\ze) / G = \{ 0 \}$

\subsection*{\nmb0{10}. Palais' treatment of $\g$-manifolds }
In \cit!{7}, Palais considered $\g$-actions on finite
dimensional manifolds $M$ in the following way. He assumed from the
beginning, that $M$ may be a non-Hausdorff manifold, since the completion
may be non-Hausdorff. Then he introduces notions which we can express as
follows in the terms introduced here:

\begin{enumerate}
\item[(\nmb.{10.1})] 
$(M,\ze)$ is called {\it generating} 
if it generates a local $G$--transformation group. 
See \cit!{7}, II,2, Def.\ V and II,7, Thm.\ XI.
This holds if and only if the leaves of the graph 
foliation on $G\x M$ described in section \nmb!{3} are Hausdorff. 
For Hausdorff $\g$--manifolds this is always the case. 
\item[(\nmb.{10.2})] 
$(M,\ze)$ is called {\it uniform} if
$\pr_1:L(e,x)\to G$ in \thetag{\nmb!{3.2}} is a covering map for each $x\in M$. 
See \cit!{7}, III,6, Def.\ VIII and III,6, Thm.\ XVII, Cor., Cor.2.
In the Hausdorff case the $\g$--action is then complete and it may 
be integrated to a $\widetilde G$--action, where $\widetilde G$ is a 
simply connected Lie group with Lie algebra $\g,$ 
so that ${_{\widetilde{G}}}M \cong M.$ 
\item[(\nmb.{10.3})] 
$(M,\ze)$ is called {\it univalent} if
$\pr_1:L(e,x)\to G$ in \thetag{\nmb!{3.2}} is injective for $\forall x.$    
See \cit!{7}, III,2, Def.\  VI and III,4, Thm.\ X.
\item[(\nmb.{10.4})] 
$(M,\ze)$ is called {\it globalizable} if there
exists a (non-Hausdorff) $G$-manifold $N$ which contains $M$ equivariantly
as an open submanifold. 
See \cit!{7}, III,1, Def.\ II and III,4, Thm.\ X.
This is a severe condition which is not satisfied
in examples \nmb!{4} and \nmb!{6} above. 
\end{enumerate}

Palais' main result on (non-Hausdorff) manifolds with a vector field says
that \thetag{\nmb!{10.1}}, \thetag{\nmb!{10.3}}, and
\thetag{\nmb!{10.4}} are equivalent. 
See \cit!{7}, III,7, Thm.\ XX.

On (non-Hausdorff) $\g$-manifolds his main result is that
\thetag{\nmb!{10.3}} and \thetag{\nmb!{10.4}} are equivalent. 
See \cit!{7}, III,1, Def.\ II and III,4, Thm.\ X, and also
III,2, Def.\ VI and III,4, Thm.\ X.

\subsection*{\nmb0{11}. Concluding remarks }

\therosteritem{\nmb.{11.1}}
A suitable setting for further development might be the class 
of {\it discrete} $\g$--manifolds, that is $\g$--manifolds for 
which the $\widetilde G$--space ${_{\widetilde G}}M$ is $T_1$, 
or equivalently the leaves of the graph foliation ${\mathcal F}_\ze$ 
on $\widetilde{G} \x M$ are closed. In this case, the charts 
$j_g:M \to {_{\widetilde G}}M$ in \thetag{\nmb!{5.1}} are 
local diffeomorphisms with respect to the unique smooth 
structure on ${_{\widetilde G}}M$ and ${_{\widetilde G}}M$ is 
a smooth manifold, albeit not necessarily Hausdorff. 

\therosteritem{\nmb.{11.2}}
In the context of \thetag{\nmb!{11.1}}, there are several 
definitions of {\it proper} $\g$--actions, all of which are 
equivalent to saying that the $\widetilde{G}$--action on 
${_{\widetilde G}}M$ is proper. Many properties of proper 
actions will carry over to this case.

\bibliographystyle{amsplain}
\providecommand{\bysame}{\leavevmode\hbox to3em{\hrulefill}\thinspace}
\providecommand{\MR}{\relax\ifhmode\unskip\space\fi MR }
\providecommand{\MRhref}[2]{%
  \href{http://www.ams.org/mathscinet-getitem?mr=#1}{#2}
}
\providecommand{\href}[2]{#2}

\end{document}